\documentclass[conference]{IEEEtran}
\IEEEoverridecommandlockouts
\usepackage{scrextend}
\usepackage[utf8]{inputenc}
\usepackage{bbm}
\usepackage{bm}
\usepackage{xcolor}
\usepackage{multirow}
\usepackage{bbold}
\usepackage{siunitx}
\usepackage{mathtools}
\usepackage{enumitem}
\usepackage{algpseudocode}
\usepackage{mathrsfs}
\usepackage[linesnumbered,ruled,vlined]{algorithm2e}

\renewcommand{\baselinestretch}{0.966}

\usepackage{cite}

\let\Right\right
\let\Left\left
\makeatletter
\def\right#1{\Right#1\@ifnextchar){\!\right}{}}
\def\left#1{\Left#1\@ifnextchar({\!\left}{}}
\makeatother

\ifCLASSOPTIONcompsoc
    \usepackage[caption=false, font=normalsize, labelfont=sf, textfont=sf]{subfig}
\else
\usepackage[caption=false, font=footnotesize]{subfig}
\fi

\ifCLASSINFOpdf
   \usepackage[pdftex]{graphicx}
\else
   \usepackage[dvips]{graphicx}
\fi

\usepackage{multirow}

\begin{document}
\title{Sample-Based Piecewise Linear Power Flow \\ Approximations Using Second-Order Sensitivities
\thanks{This work is supported by the advanced grid modeling (AGM) program of the Department of Energy, Office of Electricity and NSF award \#2145564.}}

\author{
    \IEEEauthorblockN{Paprapee Buason, Sidhant Misra}
    \IEEEauthorblockA{
        \textit{Los Alamos National Laboratory} \\
        Los Alamos, NM, USA \\
        \{buason, sidhant\}@lanl.gov
    }
    \and
    \IEEEauthorblockN{Daniel K. Molzahn}
    \IEEEauthorblockA{
        \textit{School of Electrical and Computer Engineering} \\
        \textit{Georgia Institute of Technology} \\
        Atlanta, GA, USA \\
        molzahn@gatech.edu
    }
}

\maketitle

\begin{abstract}
The inherent nonlinearity of the power flow equations poses significant challenges in accurately modeling power systems, particularly when employing linearized approximations. Although power flow linearizations provide computational efficiency, they can fail to fully capture nonlinear behavior across diverse operating conditions. To improve approximation accuracy, we propose conservative piecewise linear approximations (CPLA) of the power flow equations, which are designed to consistently over- or under-estimate the quantity of interest, ensuring conservative behavior in optimization. The flexibility provided by piecewise linear functions can yield improved accuracy relative to standard linear approximations. However, applying CPLA across all dimensions of the power flow equations could introduce significant computational complexity, especially for large-scale optimization problems. In this paper, we propose a strategy that selectively targets dimensions exhibiting significant nonlinearities. Using a second-order sensitivity analysis, we identify the directions where the power flow equations exhibit the most significant curvature and tailor the CPLAs to improve accuracy in these specific directions. This approach reduces the computational burden while maintaining high accuracy, making it particularly well-suited for mixed-integer programming problems involving the power flow equations.
\end{abstract}


\section{Introduction} \label{sec:introduction}

Describing the steady-state relationships between bus voltages and power injections, the power flow equations are fundamental to the operation, analysis, and planning of power systems.
The inherent nonlinearity of these equations introduces non-convexities that lead to challenging optimization problems. This complexity is even more pronounced when dealing with mixed-integer problems or uncertainty quantification for systems with stochastic behavior.

To mitigate these challenges, various linearization techniques are used to simplify power flow representations. Methods such as DC power flow for transmission systems~\cite{stott2009}, LinDistFlow for distribution networks~\cite{baran1989}, and first-order Taylor expansions of the AC power flow equations are commonly employed. These approximations offer computational tractability but are often based on limiting assumptions such as near-nominal voltage magnitudes and small voltage angle differences. Although such assumptions may hold in certain operating conditions, they can lead to inaccuracies when applied to more complex or broader system conditions. 

Recent research has increasingly focused on \textit{adaptive} linear power flow approximations, which are tailored to provide greater accuracy for a specific system and operating range of interest. Approaches in this area include optimization-based methods~\cite{misra2018optimal} and data-driven techniques~\cite{BUASON2022,buason2024adaptive, buason_cbla, CHEN2022108573, liu2018}; see the recent surveys in~\cite{10202779, jia2023tutorial1}. These methods aim to strike a balance between computational efficiency and accuracy by customizing linearization coefficients to better capture power flow nonlinearities. Unlike traditional approaches such as DC power flow, LinDistFlow, or first-order Taylor expansions—which prioritize simplicity and computational speed but often sacrifice accuracy—adaptive methods allocate additional computational resources upfront to calculate coefficients that align with specific operating conditions. This upfront effort proves especially advantageous in applications with an offline/online split, where coefficients can be precomputed offline using forecasted data and later applied in real-time computations. By improving the precision of linearized models, adaptive approximations enable efficient solutions to complex optimization problems, such as infrastructure planning~\cite{austgen2023comparisons, owen_aquino_talkington_molzahn-EVacuation_feeder}, AC unit commitment~\cite{castillo2016}, and sensor placement optimization~\cite{buason2023datadriven}, where traditional linearization methods may lead to significant inaccuracies.

Building on this foundation, sample-based conservative linear approximations (CLAs) were introduced to further enhance the usefulness of adaptive methods by incorporating conservativeness into the construction process~\cite{BUASON2022}. These approximations are designed to systematically over- or under-estimate specific quantities, such as voltage magnitudes, to ensure robust constraint enforcement in optimization problems. By embedding conservativeness during their formulation, CLAs improve the reliability of inequality constraints, such as those imposed by voltage or power flow limits, making them especially effective in applications where strict adherence to constraints is critical; see~\cite{buason2023datadriven,gupta_buason_molzahn-fairness_pv_limit,owen_aquino_talkington_molzahn-EVacuation_feeder} for relevant applications.

In this paper, we introduce the \textit{conservative piecewise linear approximation} (CPLA) as an extension of the CLA approach to better capture the inherent nonlinearity of the power flow equations. While CLAs provide linear approximations that are conservative by design, they are limited in their ability to adapt to highly nonlinear regions. CPLA addresses this limitation by introducing piecewise linear approximations, which enable a more precise representation of nonlinear behavior. 

To identify these nonlinear regions, we leverage a second-order sensitivity analysis to compute the directions in which the power flow equations exhibit the greatest curvature. The results from~\cite{buason2024adaptive} indicate that there are only a few dominant nonlinear directions. Thus, CPLA focuses on implementing piecewise linear approximations in a few highly nonlinear directions, while conservative linear approximations are applied in the remaining directions to manage computational complexity. This targeted approach balances improved accuracy in capturing nonlinearities with computational tractability.  

In summary, the main contributions of this paper are:
\begin{enumerate}[ label=(\textit{\roman*})]
    \item Introduction and formulation of a conservative piecewise linear approximation (CPLA) of the power flow equations based on a second-order sensitivity analysis.
    \item Variable reduction using continuity of the piecewise linearization to improve tractability of the CPLA, reducing the number of decision variables from exponential to polynomial in the number of highly nonlinear directions.
    \item Numerical studies on how the number of nonlinear directions and piecewise linearization breakpoints in the CPLAs impact approximation accuracy across test cases.
\end{enumerate}

The remainder of this paper is organized as follows. Section~\ref{sec:background} covers background material on the power flow equations as well as conservative linear approximations and second-order sensitivity analysis. Section~\ref{sec:CPLA} introduces the conservative piecewise linear approximation approach. Section~\ref{sec:simulation} provides the numerical results of our approach. Section~\ref{sec:future work} concludes the paper and discusses directions for future work. 


\section{Background} \label{sec:background}

This section provides background on the AC power flow equations, recently developed conservative linear approximations of these equations, and the second-order sensitivity analysis used to capture the power flow nonlinearities.

\subsection{The power flow equations} \label{sub:pf}

Let $V$ and $\theta$ denote the vectors of voltage magnitudes and angles. Let $P$ and $Q$ denote the active and reactive power injections. At each bus~$i$, the power flow equations are:
\begin{subequations}
\label{eq:power_flow}
\begin{align}
	P_i &= V_i^2 G_{ii} + \sum_{k \in \mathcal B_i} V_i V_k(G_{ik}\cos \theta_{ik} + B_{ik}\sin\theta_{ik}), \label{eq:dP} \\
	Q_i &= -V_i^2 B_{ii} + \sum_{k \in \mathcal B_i} V_i V_k(G_{ik}\sin\theta_{ik} - B_{ik}\cos\theta_{ik}), \label{eq:dQ}
\end{align}
\end{subequations}
where $\theta_{ik} := \theta_i - \theta_k$. The subscript $(\cdot)_{i}$ denotes a quantity at bus $i$, the subscript $(\cdot)_{ik}$ denotes a quantity from or connecting bus $i$ to $k$, and $\mathcal{B}_i$ denotes the set of all neighboring buses to bus $i$, inclusive of bus $i$ itself. Here, $G$ and $B$ represent the real and imaginary parts, respectively, of the admittance matrix.  

\subsection{Conservative linear approximations} \label{sub:CLA}

\begin{figure}[t]
	\centering 
	\includegraphics[trim=0.5cm 0.5cm 0.2cm 0.8cm, clip, width=1\linewidth]{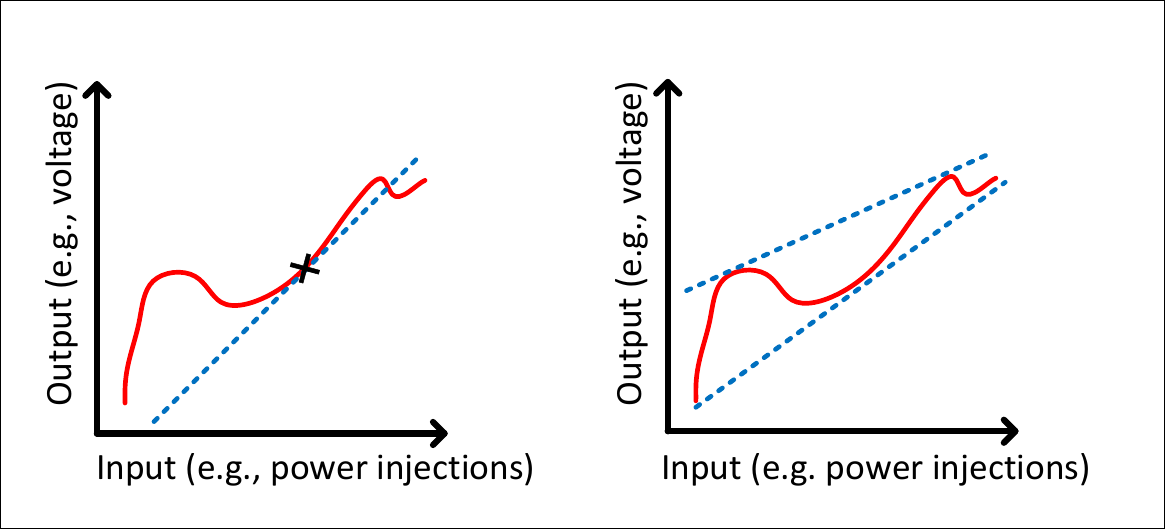} 
	\caption{The figure presents a visual comparison between a standard linear approximation (depicted on the left) and CLAs (depicted on the right). The solid line represents the nonlinear function being analyzed. In the left illustration, the dotted line portrays a conventional first-order Taylor approximation centered at point $\times$, whereas in the right illustration, the dotted line above (below) signifies an over- (under-)estimating approximation.}
	\label{fig:cla_fig}
\vspace{-1em}
\end{figure}

The nonlinear nature of power flow equations in~\eqref{eq:power_flow} presents significant challenges for optimization problems where these equations act as constraints. Conservative linear approximations (CLAs)~\cite{BUASON2022} simplify this complexity by providing linear models that over- or underestimate quantities like voltage magnitudes and current flows (Fig.~\ref{fig:cla_fig}), ensuring constraint satisfaction in optimization. By adopting a sample-based approach, CLAs can be tailored to specific power systems and operating conditions, allowing for efficient approximations within defined ranges of power injections and loads. For example, CLAs may sample loads within a predetermined range $\mathcal{S} = \{ P_{d}^{\text{min}} \leq P_{d} \leq P_{d}^{\text{max}}, \enspace Q_{d}^{\text{min}} \leq Q_{d} \leq Q_{d}^{\text{max}} \mbox{ for all } d \in \mathcal{N_D}\}$ using the uniform distribution, where $(\,\cdot\,)_{d}$ denotes the load demand, $\mathcal{N_D}$ denotes the set of all load demands, and the superscripts max (min) indicate upper (lower) limits.

In addition to simplifying the representation of the power flow equations, CLAs can implicitly incorporate the characteristics of various devices, such as transformers and smart inverters~\cite{buason2023datadriven}. This enhances CLAs accuracy in reflecting the behavior of complex systems. 

Let bold quantities denote vectors and matrices. To facilitate the integration of CLAs into optimization frameworks, we denote the quantity of interest as $\gamma$. For instance, $\gamma$ could be a particular voltage magnitude or line flow that one wishes to model as a linear function of the power injections. An overestimating CLA can be expressed as:
\begin{equation}
\gamma \leq a_{0} + \bm{a}_{1}^T\begin{bmatrix}
\bm{P} \\ \bm{Q}
\end{bmatrix},
\label{eq:cla_setup}
\end{equation}
where $a_0$ and $\bm{a}_{1}$ are coefficients of a CLA, and $T$ is transpose. Assuming that \eqref{eq:cla_setup} is indeed satisfied for all power injections $\mathbf{P}$ and $\mathbf{Q}$ in the operating range of interest, we can ensure that the constraint $\gamma \leq \gamma^{\text{max}}$ is also satisfied via enforcing a linear constraint:
\begin{equation}
a_{0} + \bm{a}_{1}^T\begin{bmatrix}
\bm{P} \\ \bm{Q}
\end{bmatrix} \leq \gamma^{\text{max}}.
\label{eq:cla_constraint}
\end{equation}

The CLA approach allows us to handle nonlinear AC power flow constraints while maintaining linear constraints within the optimization framework. To determine the coefficients $a_0$ and $\bm{a}_1$ in~\eqref{eq:cla_setup}, we formulate the following regression problem that minimizes the error between the actual quantity of interest, denoted as $\gamma$, and its linear approximation:
\begin{subequations} \label{eq:regression}
\begin{align}
&\min_{a_{0}, \ \bm{a}_{1}}  \quad \frac{1}{S}\sum_{s = 1}^S \mathcal{L}\left(\gamma_{s} - \left(a_{0} + \bm{a}_{1}^T\begin{bmatrix}
\bm{P}_s \\ \bm{Q}_s
\end{bmatrix}\right)\right) \label{eq:regression_unconstrained} \\
&\text{s.t.} \quad \gamma_{s} - \left(a_{0} + \bm{a}_{1}^T\begin{bmatrix}
\bm{P}_s \\ \bm{Q}_s
\end{bmatrix}\right) \leq 0, \quad \forall s = 1,\ldots, S.\label{eq:regression_over}
\end{align}
\end{subequations}
The subscript $(\,\cdot\,)_{s}$ represents the $s^{\text{th}}$ sample and $S$ represents the number of samples. The function $\mathcal{L}(\,\cdot\,)$ corresponds to a loss function, typically chosen to reflect the desired behavior of the approximation, such as the absolute value for $\ell_1$ loss or the square for squared-$\ell_2$ loss. The objective function in~\eqref{eq:regression_unconstrained} minimizes the mismatch between the approximated and actual quantities across all samples, while constraint~\eqref{eq:regression_over} ensures that the approximated quantity consistently overestimates the actual quantity for all samples. For underestimating CLAs, a similar process is followed, with the key distinction being the reversal of the inequality direction in~\eqref{eq:regression_over}.


\subsection{Second-order sensitivity of the power flow equations} \label{sub:sos}

Many traditional linearizations often neglect second-order and higher-order terms, assuming that the system behavior remains adequately approximated by linear models within the operating range of interest. While these linear approximations can provide reasonable estimates under certain conditions, their accuracy inherently depends on the curvature of the power flow equations within these operating ranges.

As power systems exhibit nonlinear behavior, particularly under highly varying operating conditions, overlooking second-order effects can lead to inaccuracies in the analysis results. To address this limitation, it becomes imperative to examine the second-order sensitivity of the power flow equations~\cite{buason2024adaptive}. By considering the curvature and higher-order effects, we can gain deeper insights into the system behavior and refine their approximations accordingly.

Let $g$ be a vector-valued function, denoted as $g(\bm{x})$, where $\bm{x}$ represents a vector input, yielding a vector output. To simplify the representation of variables in the power flow equations, let us define these as follows:
\begin{equation}\label{eq:x_y_def}
    \bm{x} = \begin{bmatrix}\bm{P} \\ \bm{Q} \end{bmatrix}, \quad     \bm{y} = \begin{bmatrix} \bm{\theta} \\ \bm{V}\end{bmatrix}.
\end{equation}

We rewrite the power flow equations given in~\eqref{eq:power_flow} as:
\begin{align}
\bm{x} = g(\bm{y}). \label{eq:pf_simplified}
\end{align}

Assume that a specific quantity of interest $y_k$ can be written as $f(\bm{x})$. Let $\nabla_{\bm{x}} f(\bm{x_0})$ denote the gradient of $f$ at $\bm{x} = \bm{x_0}$. The second-order Taylor approximation of $y_k$ around $\bm{x_0}$ is:
\begin{equation} \label{eq:second-order-taylor}
    y_k \approx f(\bm{x_0}) + \nabla f(\bm{x_0})^T (\bm{x} - \bm{x_0}) + \frac{1}{2} (\bm{x} - \bm{x_0})^T \bm{\Lambda}_{y_k} (\bm{x} - \bm{x_0}),
\end{equation}
where $\bm{\Lambda}_{y_k}$ is the second-order sensitivity matrix for $y_k$ around $\bm{x_0}$, which can be expressed as follows:
\begin{equation}
\label{eq:Lambda}
\bm{\Lambda}_{y_k} = 
    \begin{bmatrix}
        \cfrac{\partial^2 y_k}{\partial x_1 \partial x_1} &\cfrac{\partial^2 y_k}{\partial x_1 \partial x_2} &\cdots &\cfrac{\partial^2 y_k}{\partial x_1 \partial x_{2N}} \\
        \vdots &\vdots &\ddots &\vdots \\
        \cfrac{\partial^2 y_k}{\partial x_{2N} \partial x_1} &\cfrac{\partial^2 y_k}{\partial x_{2N} \partial x_2} &\cdots &\cfrac{\partial^2 y_k}{\partial x_{2N} \partial x_{2N}} 
    \end{bmatrix}.
\end{equation}

The explicit form of the second-order sensitivity for voltages and its derivation can be found in~\cite{buason2024adaptive}. We will use the second-order sensitivities to tailor the CLPA's piecewise linearizations based on the directions with the most significant nonlinearities.

\section{Conservative Piecewise Linear Approximations} \label{sec:CPLA}

This section introduces conservative piecewise linear approximation (CPLAs) to improve the accuracy of power flow linearizations beyond what the CLA can offer. While CPLAs can provide more accurate representations by introducing multiple linear segments, applying CPLAs along all dimensions of the power flow equations can become computationally expensive, especially for large-scale optimization problems. In this section, we present the CPLA regression formulation and an efficient methodology that focuses on critical directions and variable reduction to lower the computational complexity.

\subsection{Overall computational steps and parallel scalability} \label{sub:computation}

\begin{figure}[b]
\vspace{-0.8em}
	\centering 
	\includegraphics[trim=0.4cm 0.3cm 0.4cm 0.3cm, clip, width=0.53\linewidth]{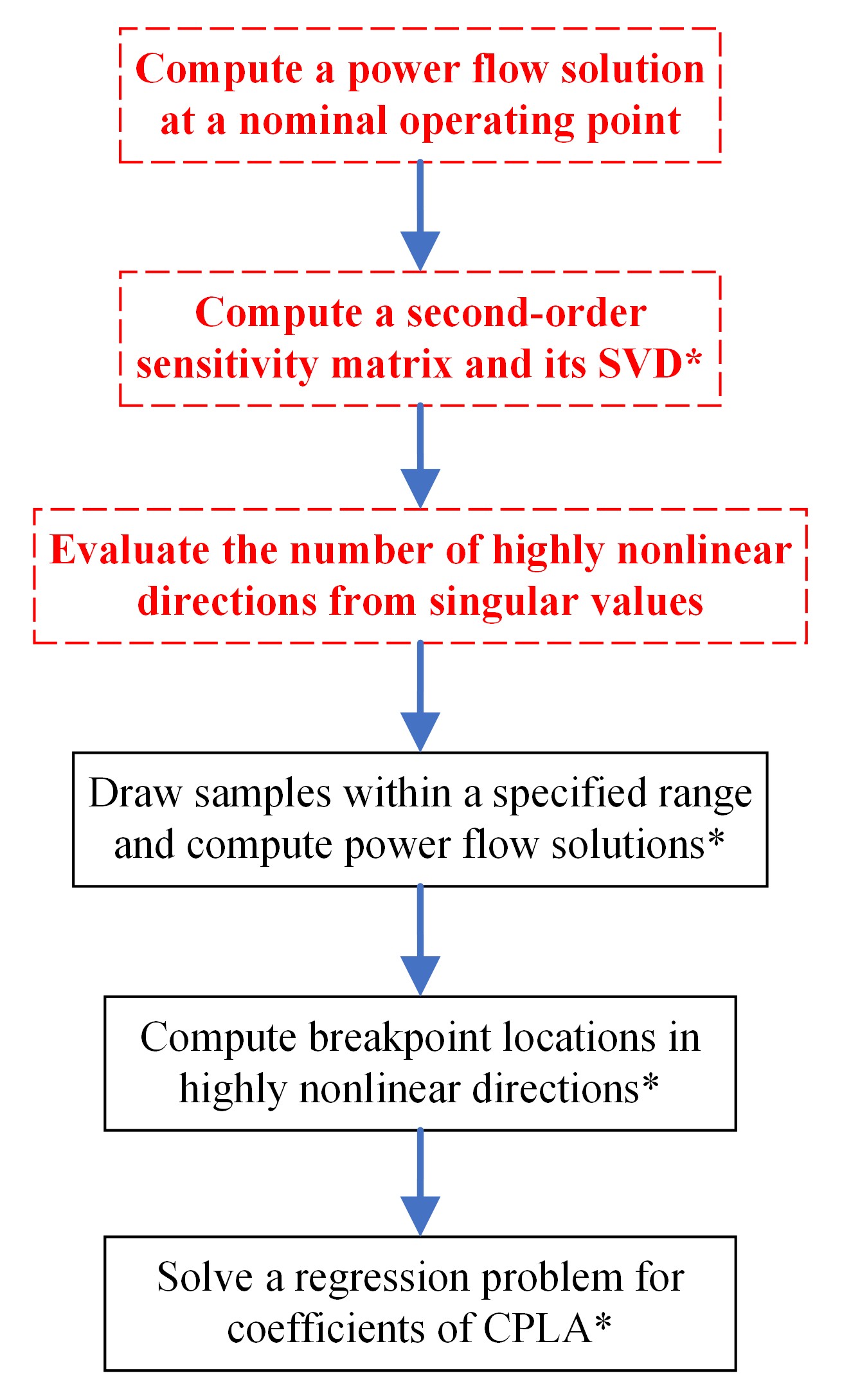} 
	\caption{Flowchart illustrating the computation processes for the CPLA method. Red-dashed boxes highlight the steps involving the second-order sensitivity matrix. Steps marked with $*$ indicate parallelizable processes.}
	\label{fig:algorithm}
\end{figure}

The CPLA computation consists of two main steps: identifying highly nonlinear directions via a second-order sensitivity analysis and solving the regression problem in~\eqref{eq:regression_CPLA} to construct piecewise linear approximations. These steps ensure that the approximations are tailored to accurately capture the nonlinear behavior of power flow equations. The flowchart in Fig.~\ref{fig:algorithm} illustrates the overall steps involved in computing the CPLA.

The approach is highly parallelizable, making it efficient for large-scale systems. The second-order sensitivity matrix computations for each bus and the power flow calculations for each sample are independent. Thus, they can be performed in parallel. Similarly, the regression problems for determining approximation coefficients are tractable linear problems that are independent across buses, enabling further parallelization. Furthermore, voltage solutions obtained from power flow computations can be reused across buses, reducing redundant calculations. This parallel structure ensures scalability and computational efficiency for large systems.

\subsection{Problem formulation} \label{sub:CPLA_formulation}

\begin{figure}[bt!]
    \vspace{-0.5em}
	\centering 
	\includegraphics[trim=0.5cm 0.2cm 0.5cm 0.4cm, clip, width=0.7\linewidth]{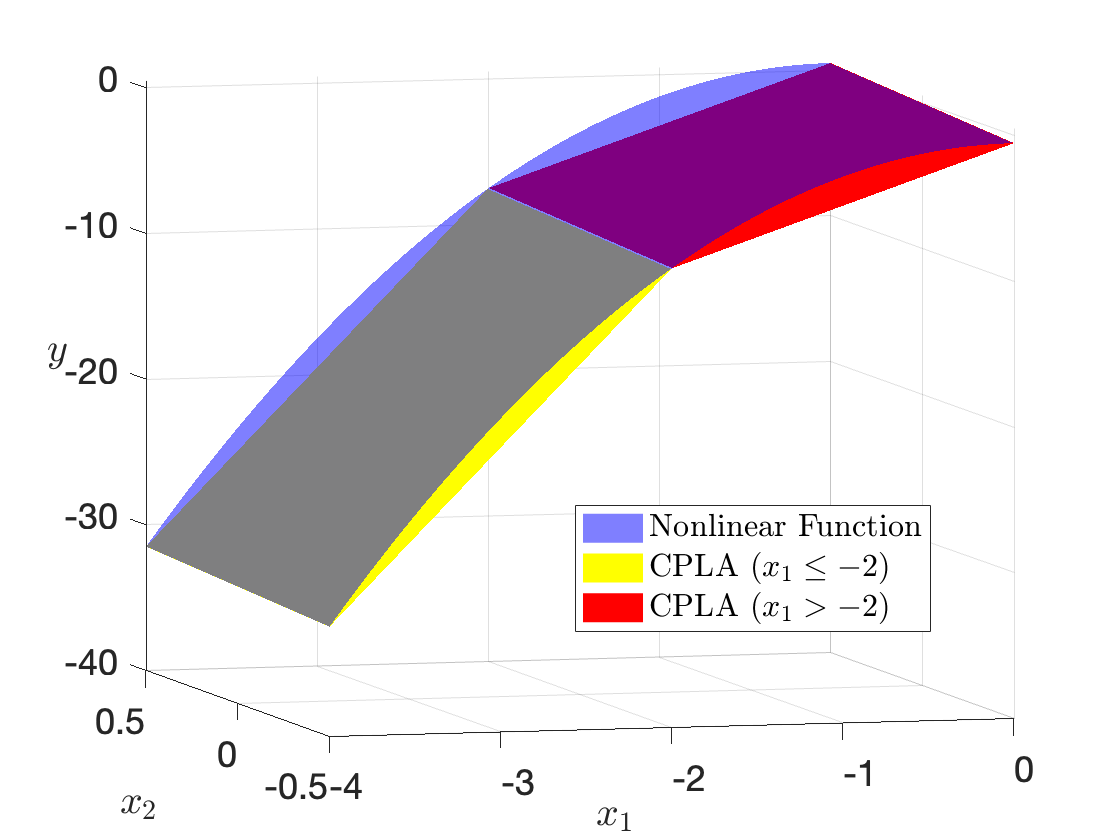} 
	\caption{An example of a CPLA is shown, where the yellow and red planes (each referred to as a ``region" in this paper) underestimate the quadratic function $y = -4x_1^2 + x_2$ (blue manifold). The breakpoint is at $x_1 = -2$, with $x_1 \leq -2$ and $x_1 > -2$ referred to as two segments. The yellow (red) plane corresponds to the function $y = 12x_1 + x_2 + 16$ ($y = 4x_1 + x_2 $).}
	\label{fig:cpla_3d}
	\vspace{-1em}
\end{figure}

CPLAs enhance the accuracy of power flow linearizations by addressing nonlinear behavior through multiple linear segments. However, not all dimensions of the power flow equations exhibit significant nonlinearity. Most directions display only mild curvature, while a select few are predominantly nonlinear and contribute substantially to the overall system behavior. To conceptually illustrate this concept, Fig.~\ref{fig:cpla_3d} presents an underestimating CPLA for a quadratic function with one highly nonlinear direction ($x_1$) and one linear direction ($x_2$). In this case, a piecewise linear function is only necessary along $x_1$, as $x_2$ exhibits linear curvature, highlighting the efficiency of focusing on critical dimensions.

Recognizing which dimensions require piecewise linearization is crucial for improving computational efficiency without sacrificing accuracy. This is where second-order sensitivity analysis becomes invaluable. By examining the second-order sensitivities of the power flow equations, we can identify and prioritize directions with pronounced nonlinear effects.

Let $[\bm{P'}_s; \bm{Q'}_s]$ represent the vector of rotated power injections (discussed in Section~\ref{sub:rotating_svd}) corresponding to highly nonlinear directions calculated from a second-order sensitivity analysis.
Similar to the regression problem in~\eqref{eq:regression}, the formulation for computing a sample-based overestimating CPLA is as follows:
\begin{subequations} \label{eq:regression_CPLA}
\begin{align}
&\min_{h(\,\cdot\,)} \quad \ \frac{1}{S}\sum_{s = 1}^S \mathcal{L}\left(\gamma_{s} - h\left(\begin{bmatrix}\bm{P}_s \\ \bm{Q}_s\end{bmatrix}, \begin{bmatrix}\bm{P'}_s \\ \bm{Q'}_s\end{bmatrix}, b_s\right)\right) \label{eq:regression_unconstrained_CPLA} \\
&\text{s.t.} \;
\gamma_{s} - h\left(\begin{bmatrix}\bm{P}_s \\ \bm{Q}_s\end{bmatrix}, \begin{bmatrix}\bm{P'}_s \\ \bm{Q'}_s\end{bmatrix}, b_s\right) \leq 0, \label{eq:regression_over_CPLA} \\
&\hphantom{\text{s.t.} \;}  h\left(\begin{bmatrix}\bm{P}_s \\ \bm{Q}_s\end{bmatrix}, \begin{bmatrix}\bm{P'}_s \\ \bm{Q'}_s\end{bmatrix}, b_s\right) = a_{0} + \bm{a}_{1}^T\begin{bmatrix}
\bm{P}_s \\ \bm{Q}_s
\end{bmatrix} + g\left(\begin{bmatrix}\bm{P'}_s \\ \bm{Q'}_s\end{bmatrix}, b_s\right), \label{eq:equality_CPLA} \\
&\hphantom{\text{s.t.} \;} \forall s = 1,\ldots, S, \notag
\end{align}
\end{subequations}
where $b_s$ denotes the segment which the $s^{\text{th}}$ sample belongs to, $h(\,\cdot\,)$ represents the CPLA that consists of a CLA and a purely CPLA, denoted as $g(\,\cdot\,)$, which represents the piecewise linear function for which we optimize the coefficients across all segments (shown in Section~\ref{sub:continuity}). Note that breakpoint locations in each direction are determined by evenly partitioning the sample space, allowing samples to share a segment in one direction but belong to different segments in another. By determining the breakpoint locations based on the drawn samples before solving for a CPLA in~\eqref{eq:regression_CPLA}, we can effectively assign each sample to its corresponding segment. 

We need to ensure that there are multiple samples within each region of the piecewise linear approximation.
To avoid overfitting in the piecewise linear approximation, it is crucial to ensure that sufficient sample data exists within each region defined by the multidimensional piecewise linear functions. Overfitting is particularly likely when multiple highly nonlinear directions and several breakpoints in each direction are considered, as regions with sparse or no samples can lead to poorly generalized approximations. To address this, we propose the following two approaches:
\begin{enumerate} 
    \item Ensure an adequate number of drawn samples. Let $C$ represent the number of regions and $S$ represent the number of samples. Define $\epsilon$ as a predefined small positive confidence value. We ensure that the probability of the undesirable event where there are one or more regions with no samples is bounded by $\epsilon$:
    \begin{align} \label{eq:bound_1} 
    C\left(1 - \cfrac{1}{C}\right)^S \leq \epsilon. 
    \end{align}
    
    Assuming $C$ is large, the bound on $S$ is:
    \begin{align} \label{eq:bound} 
    S \geq C \log \cfrac{C}{\epsilon}. 
    \end{align}

    Note that this inequality approximates the probability that there is a region with no sample. This approximation is used to simplify the computation of the inequality in~\eqref{eq:bound_1}.

    \item Use regularization by enforcing convexity or concavity in each direction, which can be determined by analyzing the singular values obtained from SVD. This approach reduces the required samples and, consequently, the computation time compared to the first approach.
\end{enumerate}

\subsection{Variable reduction from continuity} \label{sub:continuity}

In this section, we consider a piecewise linear function in $N$ dimensions with $M$ breakpoints, resulting in $M + 1$ segments within each dimension. Let $t_1, t_2, \ldots, t_N$ represent the $N$ directions and $b_1, b_2, \ldots, b_N$ with $b_i \in \{0, 1, 2, \ldots, M\}$ represent the index of the segment in each direction. We define $\bm{t} = [t_1, t_2, \ldots, t_N]^T$ and $\bm{b} = [b_1, b_2, \ldots, b_N]$. The piecewise linear function $g(\bm{t},\bm{b})$ is expressed as follows:
\begin{align}
g(\bm{t},\bm{b}) = & a_0^{b_1 b_2 \ldots b_N} + a_1^{b_1 b_2 \ldots b_N} t_1 
+ \ldots + a_N^{b_1 b_2 \ldots b_N} t_N, \label{eq:piecewise_function}
\end{align}
where the coefficients $a_{0}^{b_1 b_2 \ldots b_N}, a_1^{b_1 b_2 \ldots b_N}, \ldots, a_N^{b_1 b_2 \ldots b_N}$ correspond to $a_0, a_1, \ldots, a_N$ for the region defined by the segments $b_1, b_2, \ldots, b_N$. For instance, in Fig.~\ref{fig:cpla_3d}, since there is only one nonlinear direction $x_1$, we have a single segment index $b_1 \in \{0,1\}$, representing two segments. Specifically, $a_0^0 = 16$ and $a_1^0 = 12$ for the first segment ($b_1 = 0$), which corresponds to $x < -2$. The values of $t_1, t_2, \ldots, t_{N}$ fall within the specified range defined by the breakpoints.

We begin by examining the first breakpoint in the $t_1$ direction, denoted as $t_1^0$. At the point where $t_1 = t_1^0$, the continuity condition can be stated as follows:
\begin{align}
    a_0^{0 b_2 \ldots b_N} + a_1^{0 b_2 \ldots b_N} t_1^0 + \ldots & + a_N^{0 b_2 \ldots b_N} t_N = \notag \\
    a_0^{1 b_2 \ldots b_N} + a_1^{1 b_2 \ldots b_N} t_1^0 + \ldots & + a_N^{1 b_2 \ldots b_N} t_N.
\end{align}

Therefore, we can conclude that:
\begin{subequations}
\begin{align}
    a_0^{0 b_2 \ldots b_N} + a_1^{0 b_2 \ldots b_N} t_1^0 &= a_0^{1 b_2 \ldots b_N} + a_1^{1 b_2 \ldots b_N} t_1^0, \label{eq:a0_0}\\
    a_2^{0 b_2 \ldots b_N} &= a_2^{1 b_2 \ldots b_N}, \label{eq:a2_0}\\
    &\vdots \notag \\
    a_{N}^{0 b_2 \ldots b_N} &= a_{N}^{1 b_2 \ldots b_N} \label{eq:a(N)_0}.
\end{align}
\end{subequations}
We apply a similar approach to all other breakpoints in the $t_1$ direction. For instance, considering the $m^{\text{th}}$ breakpoint yields:
\begin{subequations}
\begin{align}
    a_0^{(m - 1) b_2 \ldots b_N} +  a_1^{(m - 1) b_2 \ldots b_N} t_1^{m - 1} &= a_0^{m b_2 \ldots b_N} +  a_1^{m b_2 \ldots b_N} t_1^{m - 1}, \label{eq:a0_(m-1)}\\
    a_2^{(m - 1) b_2 \ldots b_N} &= a_2^{m b_2 \ldots b_N}, \label{eq:a2_(m-1)}\\
    &\vdots \notag \\
    a_{N}^{(m - 1) b_2 \ldots b_N} &= a_{N}^{m b_2 \ldots b_N} \label{eq:a(N)_(m-1)}.
\end{align}
\end{subequations}
By combining \eqref{eq:a2_0}--\eqref{eq:a(N)_0} and \eqref{eq:a2_(m-1)}--\eqref{eq:a(N)_(m-1)} for all $m = 0,1,2,\ldots,M - 1$, we obtain:
\begin{subequations}
\begin{align}
    a_0^{b_1' b_2 \ldots b_N} &= \sum_{i = 0}^{b_1' - 1} t_1^i (a_0^{i b_2 \ldots b_N} - a_0^{(i + 1) b_2 \ldots b_N}) + a_0^{0 b_2 \ldots b_N}, \label{eq:a0_i1}\\
    a_2^{b_1 b_2 \ldots b_N} &= a_2^{b_1' b_2 \ldots b_N}, \label{eq:a2_i1}\\
    &\vdots \notag \\
    a_{N}^{b_1 b_2 \ldots b_N} &= a_{N}^{b_1' b_2 \ldots b_N} \label{eq:a(N)_i1},
\end{align}
\end{subequations}
for any $b_1 \in \{0,1,2,\ldots,M\}$ and $b_1' \in \{ 1,2,\ldots,M\}$.

We first focus specifically on the coefficient $a_2$ and make a claim regarding the other coefficients. By examining the breakpoints in the direction of $t_3$ and following an approach similar to the one outlined in equations~\eqref{eq:a2_i1} to~\eqref{eq:a(N)_i1}, we derive an equality for the coefficient $a_2$ as follows:
\begin{align}
    a_2^{b_1 b_2 b_3 \ldots b_N} &= a_2^{b_1 b_2 b_3' \ldots b_N}, \label{eq:a2_i3}
\end{align}
for any $b_3, b_3' \in \{0,1,2,\ldots,M\}$.

Combining \eqref{eq:a2_i1} with \eqref{eq:a2_i3}, we obtain:
\begin{align}
    a_2^{b_1 b_2 b_3 \ldots b_N} = a_2^{b_1' b_2 b_3 \ldots b_N} =  a_2^{b_1 b_2 b_3' \ldots b_N}. \label{eq:a2_all}
\end{align}

Proceeding with the analysis, if we perform the same simplification across all directions except for direction $t_2$, we find that the coefficient $a_2$ is solely dependent on the segment $b_2$. Similarly, each coefficient $a_n$ depends exclusively on its respective segment $b_n$ for all $n = 1, 2, \ldots, N$. 

Now, we analyze the coefficient $a_0$ and observe how it changes with directions. Utilizing the fact that $a_n$ is dependent solely on direction $t_n$, coupled with information from~\eqref{eq:a0_i1}, $a_0$ is given by:
\begin{equation} \label{eq:a0_all}
a_0^{b_1 b_2 \ldots b_N} = \sum_{k = 1}^{N} \left(\sum_{b_k' = 0}^{b_k - 1} t_k^{b_k'}(a_k^{b_k'} - a_k^{b_k' + 1})\right) + a_0^{00 \ldots 0}.
\end{equation}

Consequently, the piecewise linear function in~\eqref{eq:piecewise_function} can be represented in a simplified form as:
\begin{align} \label{eq:complete_piecewise_function}
g(\bm{t},\bm{b}) = a_0^{b_1 b_2 \ldots b_N} + a_1^{b_1} t_1 + a_2^{b_2} t_2 + \ldots + a_{N}^{b_N} t_N.
\end{align}

Leveraging continuity allows independent selection of breakpoints and simplifies the decision variables from $(N+1)\times(M+1)^{N}$ to $(N+1)\times(M+1)$.

\subsection{Rotating coordinates from singular value decomposition} \label{sub:rotating_svd}

The directions of high curvature are identified using the significant singular values obtained from singular value decomposition (SVD) of the second-order sensitivity matrix $\mathbf{\Lambda}$ in~\eqref{eq:Lambda}. The singular vectors corresponding to these significant singular values are the directions in which we implement the CPLA.
Let $\bm{\Lambda} = \bm{U}\bm{S}\bm{U}^T$ be the singular value decomposition of $\bm{\Lambda}$ and $\widetilde{\bm{U}}$ be the sub-matrix of $\bm{U}$ obtained by selecting the columns corresponding to the significant singular values. Let $\bm{w}$ denote the complete vector of power injections. The vector of rotated power injections, $\bm{t}$ in~\eqref{eq:complete_piecewise_function}, which represents the directions of high curvature is given by $\bm{t} = \widetilde{\bm{U}}^T\bm{w}$.

Since the segment indices $b_i$ are automatically determined from the samples, the regression problem in~\eqref{eq:regression_unconstrained_CPLA} using the representation of $g(\,\cdot\,)$ in~\eqref{eq:complete_piecewise_function} has a polyhedral constraint set and can be solved efficiently for several loss functions $\mathcal{L}(\,\cdot\,)$.


\section{Numerical results} \label{sec:simulation}

In this section, we present numerical results to evaluate the effectiveness of the CPLA approach for approximating nonlinear power flow equations. The analysis is structured into three parts: examining the impact of singular vectors associated with dominant singular values, analyzing the effect of varying the number of segments, and exploring the combined influence of these factors to optimize performance.

Simulations are conducted on the \textit{case30}, \textit{case33bw}, \textit{case141}, and \textit{case2383wp} test cases from M{\sc atpower}~\cite{zimmerman_matpower_2011}, with optimization problems solved using the YALMIP toolbox~\cite{Lofberg_Yalmip}. For each case, 10,000 samples are drawn by varying power injections within 50\% to 150\% of their nominal values. Reported voltages are in per unit (pu), and the $\ell_1$ norm, denoted by $\mathcal{L}(\,\cdot\,)$, is used as the loss function. All computations are performed on a MacBook Pro with an Apple M1 Pro chip (10 cores, 16~GB RAM).

\subsection{Effects of nonlinear directions and breakpoints} \label{sub:sim_direction}

Piecewise linear approximation accuracy is strongly influenced by the specific nonlinear directions being considered. These directions are derived from the singular vectors corresponding to the dominant singular values of the system, with each singular vector representing a unique nonlinear direction. The amount of curvature in a given direction, as captured by its corresponding singular value, determines how well the approximation performs along that direction. Highly nonlinear directions are more challenging to approximate accurately. Analyzing  one direction at a time provides insights into their influence on approximation quality.

In addition to the choice of nonlinear directions, the number of breakpoints used to represent a nonlinear function is another critical factor influencing the accuracy of a piecewise linear approximation. A finer segmentation can capture the underlying nonlinearity more effectively, resulting in significantly reduced approximation errors. However, this comes with the trade-off of increased computational demands. In this section, we analyze how the individual nonlinear directions and the number of breakpoints impact approximation accuracy.

\begin{figure}[b!]
\vspace{-2em}
    \centering 
    \subfloat[One singular vector\label{fig:30bus_eigen}]{\includegraphics[trim=0.35cm 0cm 0.4cm 0.1cm, clip, width=0.5\linewidth]{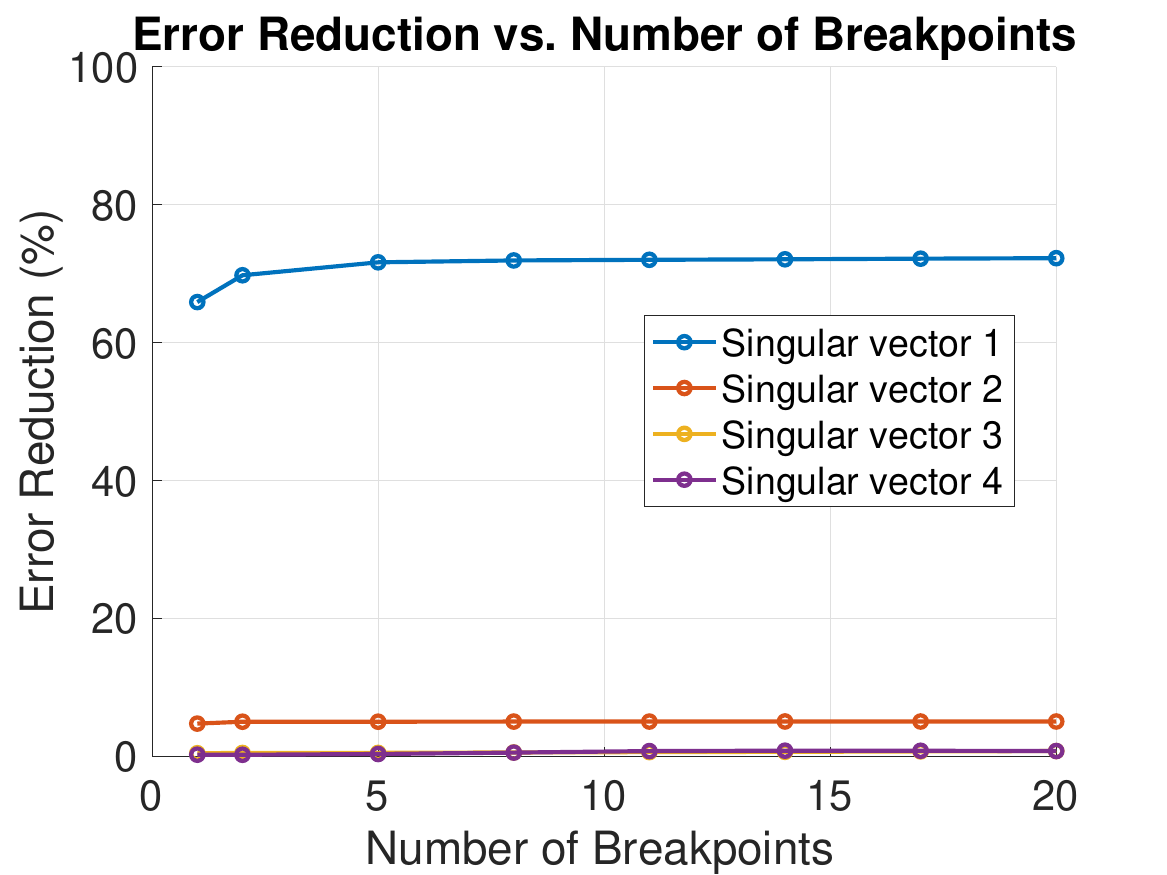}}
    \hfill
    \subfloat[Two singular vectors\label{fig:141bus_eigen}]{\includegraphics[trim=0.35cm 0cm 0.4cm 0.1cm, clip, width=0.5\linewidth]{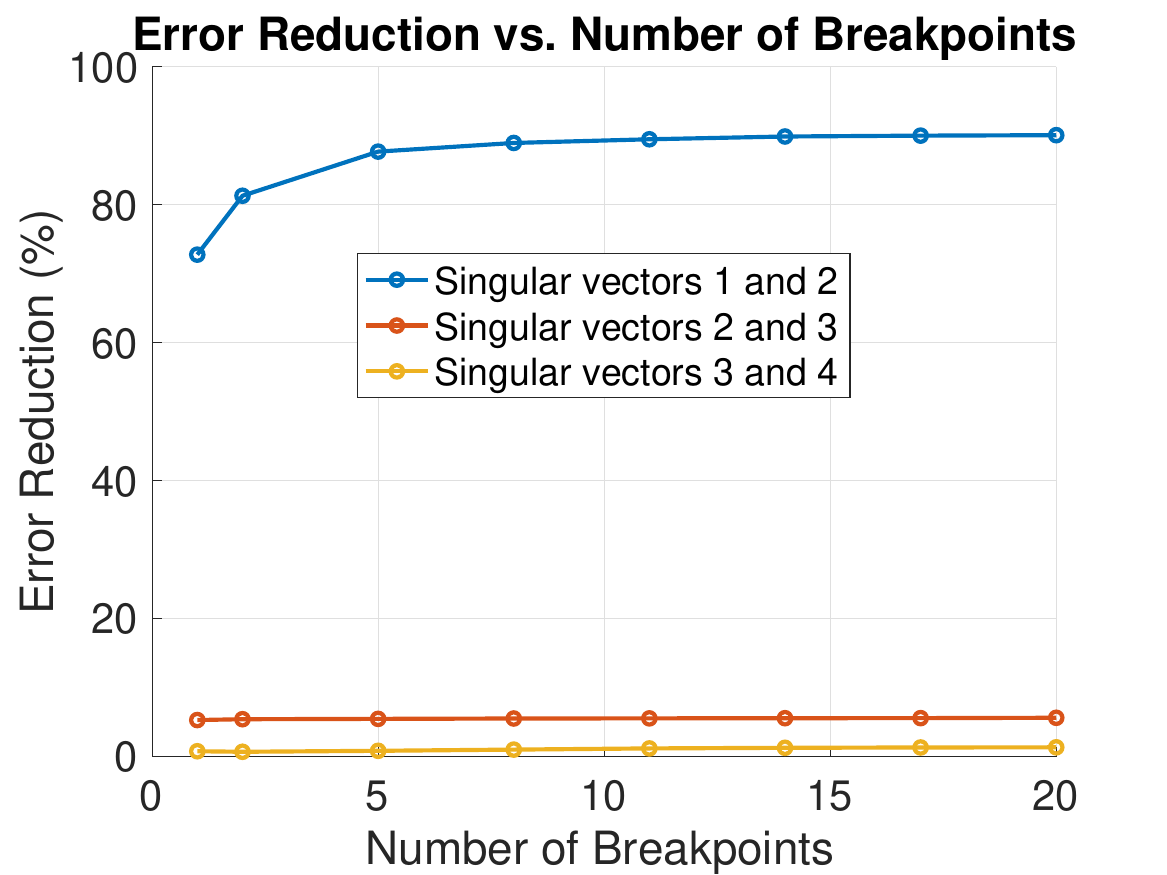}}
    \caption{Percentage error reduction of voltage magnitudes vs. number of breakpoints in the 141-bus system at bus 80 when considered (a) one singular vector and (b) two singular vectors.}
    \label{fig:error_reduction_comparison}
\end{figure}

\begin{table*}[t!] 
\caption{Approximation errors for underestimating voltage magnitudes at a specific bus}
\label{table:voltage_CPLA}
\centering
\setlength\tabcolsep{2.5pt}
\begin{tabular}{c|c|c|c|c|c|c|c|c|c|c}
  \multirow{3}{*}{Cases} & \multirow{3}{*}{Bus} & \multirow{3}{*}{$t_{\text{pf}}$ [ms]} & \multicolumn{2}{c|}{CLA} & \multicolumn{6}{c}{CPLA} \\ 
  \cline{4-11}
  & & & Error/sample & time & \multicolumn{3}{c|}{Error/sample for each \#breakpoint(s)$^*$ [pu]} & \multicolumn{3}{c}{$t_{\text{CPLA}}$ [s]} \\
  \cline{6-11}
  & & & [pu] & [s] & 1 & 5 & 10 & 1 & 5 & 10  \\
  \hline \hline
    \multirow{2}{*}{\textit{case30}} & \multirow{2}{*}{20} & \multirow{2}{*}{2.26}  & \multirow{2}{*}{$7.36\times10^{-3}$} & \multirow{2}{*}{1.51} & $4.32\times10^{-3}$ & $4.24\times10^{-3}$  & $4.21\times10^{-3}$ & \multirow{2}{*}{1.67} & \multirow{2}{*}{1.72} & \multirow{2}{*}{1.73} \\
    & & & &  & (41.30\%) & (42.42\%) & (42.84\%) & & & \\
  \hline
  \multirow{2}{*}{\textit{case33bw}} & \multirow{2}{*}{33} &  \multirow{2}{*}{2.43}  & \multirow{2}{*}{$2.84\times10^{-4}$} & \multirow{2}{*}{2.27} & $1.37\times10^{-4}$ & $1.18\times10^{-4}$  & $1.16\times10^{-4}$ & \multirow{2}{*}{2.37} & \multirow{2}{*}{2.38} & \multirow{2}{*}{2.41}  \\
  & & & &  & (51.61\%) & (58.51\%) & (59.13\%) & & & \\
  \hline
  \multirow{2}{*}{\textit{case141}} &  \multirow{2}{*}{80} &    \multirow{2}{*}{5.13}  & \multirow{2}{*}{$6.75\times10^{-5}$} & \multirow{2}{*}{5.57} & $2.31\times10^{-5}$ & $1.92\times10^{-5}$  & $1.90\times10^{-5}$ & \multirow{2}{*}{5.95} & \multirow{2}{*}{7.98} & \multirow{2}{*}{8.16} \\
  & & & &  & (65.86\%) & (71.63\%) & (71.93\%) & & & \\
  \hline
   \multirow{2}{*}{\textit{case2383wp}} &  \multirow{2}{*}{466} &    \multirow{2}{*}{169.52} &  \multirow{2}{*}{$7.45\times10^{-5}$} &  \multirow{2}{*}{125.16}  &  $4.31\times10^{-5}$ & $3.64\times10^{-5}$  & $3.59\times10^{-5}$ & \multirow{2}{*}{148.24} & \multirow{2}{*}{175.86} & \multirow{2}{*}{182.47}  \\
    & & & &  & (41.83\%) & (50.81\%) & (51.43\%) & & & \\
  \hline
\end{tabular}\\[0.2em]
\footnotesize{$t_{\text{pf}}$ = Average time required to solve 1 power flow solution in milliseconds. \\$t_{\text{CPLA}}$ = Time required to compute the second-order sensitivity matrix, perform the SVD, and compute the CPLA for a specified number of breakpoints.
\\ *The percentage reduction in errors compared to the errors from conservative linear approximation (CLA).}
\vspace{-1.5em}
\end{table*}

Fig.~\ref{fig:error_reduction_comparison} shows the error reduction compared to the CLA results as a function of the number of breakpoints for the 141-bus system at bus 80. Different lines represent the results obtained using various singular vectors. The singular vector associated with the most dominant singular value yields the greatest accuracy improvement. Increasing the number of breakpoints significantly reduces approximation errors for all singular vectors; however, this improvement diminishes beyond five breakpoints, where the error reduction stabilizes. For the most dominant singular vector, the error reduction reaches 72.22\%, highlighting its substantial contribution to improving CPLA accuracy (refer to Fig.~\ref{fig:30bus_eigen}). The improvement is particularly pronounced for the singular vector associated with the most dominant singular value, which achieves a significant error reduction. In contrast, other singular vectors provide smaller improvements, with the second singular vector achieving a maximum error reduction of about 5\%, while subsequent singular vectors contribute even less.

Considering two singular vectors, Fig.~\ref{fig:141bus_eigen} illustrates the error reduction for the CPLA of voltages at bus 80 of the 141-bus system. The results show that implementing CPLA in the directions of the first two singular vectors, corresponding to the two most dominant singular values, achieves a combined improvement of up to 90.07\%. In particular, incorporating two singular vectors provides additional improvement over using only the first singular vector, demonstrating the advantage of capturing multiple high-curvature directions.

\begin{figure}[t]
    \vspace{-0.5em}
	\centering 
	\includegraphics[trim=0.1cm 0.1cm 0.1cm 0.1cm, clip, width=0.58\linewidth]{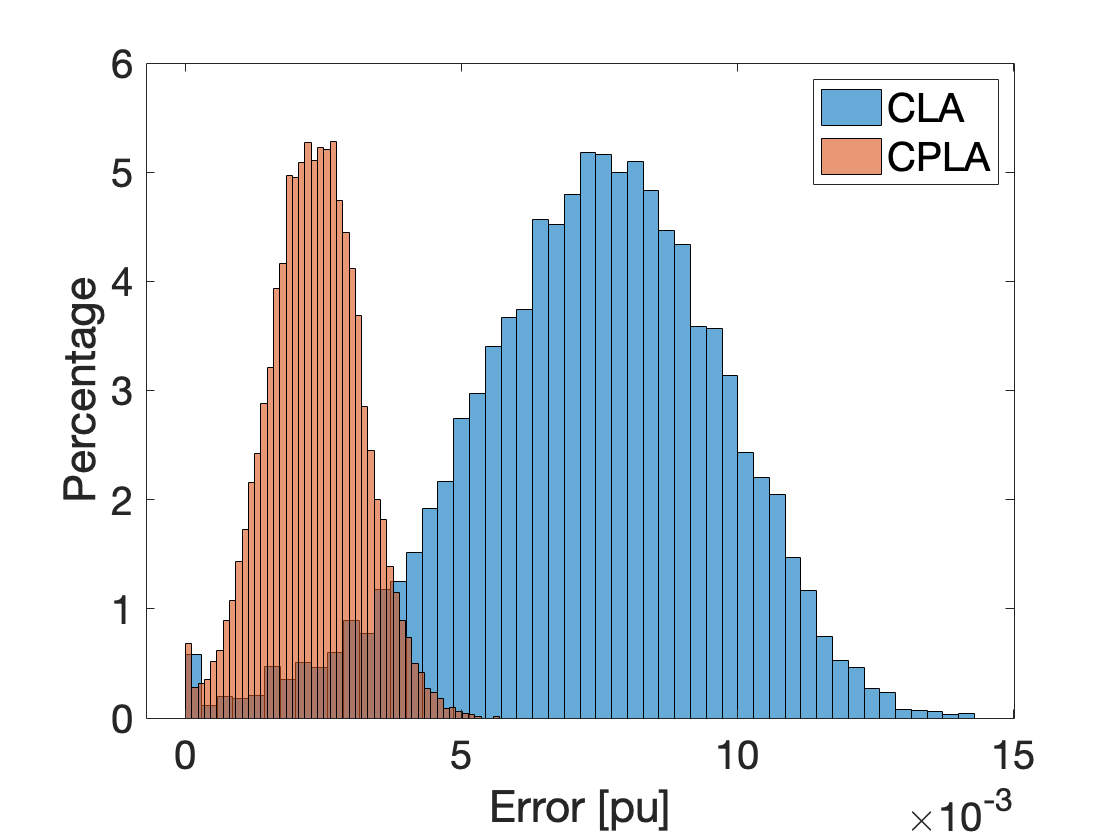} 
	\caption{Error histograms from the underestimating CLA (blue) and the underestimating CPLA (red) for voltage magnitudes at bus 20 in the IEEE 30-bus system.}
	\label{fig:cpla_cla_case30}
	\vspace{-1.2em}
\end{figure}

Table~\ref{table:voltage_CPLA} presents the detailed approximation errors, error reductions, and computation times for computing the underestimating voltages by selecting the singular vector corresponding to the most dominant singular value in the \textit{case30}, \textit{case33bw}, \textit{case141}, and \textit{case2383wp} test cases at a specific bus. The results demonstrate that the proposed piecewise linear approach reduces errors by 41.30\% to 71.93\% while requiring comparable computation time to CLA, highlighting its effectiveness in balancing improved accuracy with computational efficiency. For better visualization, the plot comparing the results of CLA and CPLA with five breakpoints is shown in Fig.~\ref{fig:cpla_cla_case30}, demonstrating that CPLA provides a better approximation than CLA.

\subsection{Singular values and their effects to the nonlinear directions} \label{sub:sim_multiple_direction}

The insights from analyzing nonlinear directions and segmentation granularity from the previous section highlight the trade-offs in achieving accurate piecewise linear approximations. Significant accuracy gains come from dominant nonlinear directions, with diminishing returns from subsequent singular vectors. Similarly, increasing breakpoints improves accuracy only up to a point, beyond which additional segments offer minimal benefit. These findings emphasize the need to balance dominant nonlinear directions with segmentation granularity. This section explores the impact of varying the number of nonlinear directions while keeping the number of breakpoints fixed and examines the trend of singular values.

\begin{figure}[t]
\vspace{-1em}
    \centering 
    \subfloat[Error reduction\label{fig:141bus_5breakpoints}]{\includegraphics[trim=0.4cm 0cm 0.4cm 0.1cm, clip, width=0.5\linewidth]{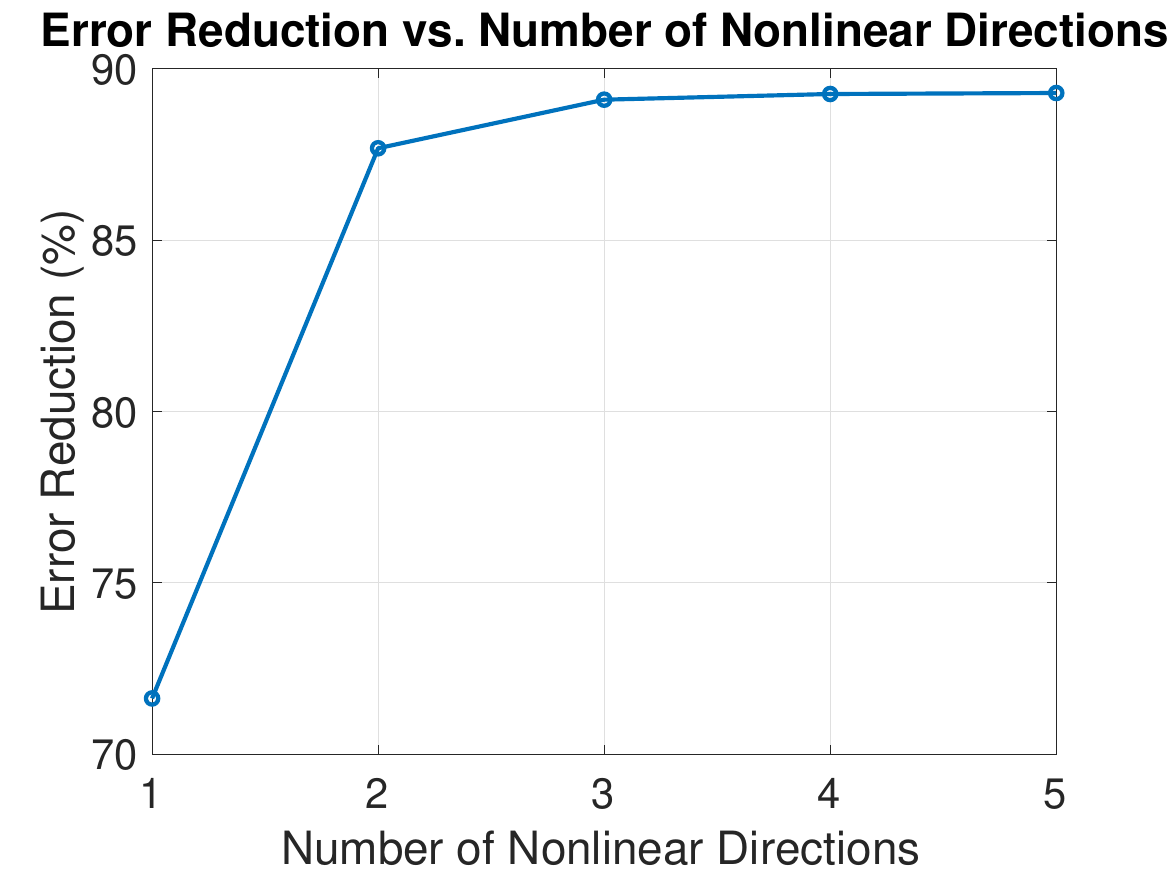}}
    \hfill
    \subfloat[Singular values\label{fig:case141_singular_values}]{\includegraphics[trim=0.4cm 0cm 0.4cm 0.1cm, clip, width=0.5\linewidth]{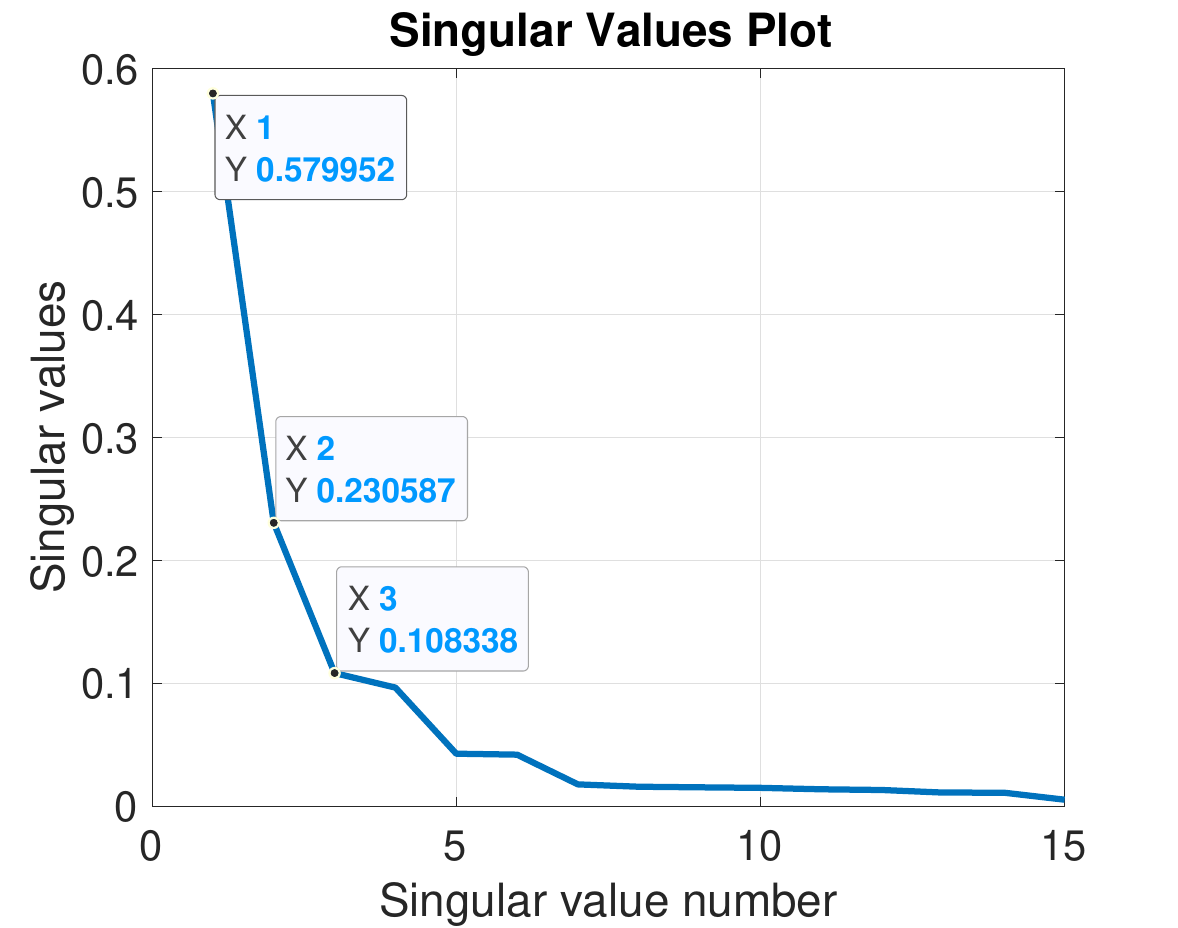}}
    \caption{(a) Percentage error reduction of voltage magnitudes vs. the number of nonlinear directions and (b) plot of the first 15 most dominant singular values for the 141-bus system at bus 80.}
    \label{fig:case141_multiple_sv}
    \vspace{-1em}
\end{figure}

Fig.~\ref{fig:case141_multiple_sv} illustrates the reduction in error when varying the number of nonlinear directions with five breakpoints in each direction, alongside the first 15 most dominant singular values. The results in Fig.~\ref{fig:141bus_5breakpoints} indicate that incorporating up to five nonlinear directions achieves a significant error reduction of 89.29\%. This performance is comparable to the results achieved by considering two nonlinear directions with ten breakpoints in each direction, as shown in Fig.~\ref{fig:141bus_eigen}. The singular values plotted in Fig.~\ref{fig:case141_singular_values} reveal a steep decline, dropping from 0.58 for the most dominant direction to 0.11 for the third, and continuing to decrease for subsequent directions. This shows that only a few highly nonlinear directions contribute significantly to reducing errors, emphasizing the efficiency of focusing computations on these dominant directions.

\section{Conclusion and future work} \label{sec:future work}

This paper demonstrates the effectiveness of the CPLA approach for approximating nonlinear power flow equations, emphasizing a balance between accuracy and computational efficiency. The results highlight that only the first few singular vectors, corresponding to the most dominant singular values, contribute significantly to improving the accuracy of the approximation. This observation not only reinforces the importance of targeting key nonlinear directions but also reduces computational demands by limiting the focus to a small subset of directions. Furthermore, while increasing the number of segments improves precision, the diminishing returns observed beyond a certain point emphasize the need for an optimal segmentation strategy. This approach provides a practical and efficient solution for addressing the challenges of nonlinear power flow equations in optimization problems.

In future work, we aim to apply our CPLA method to power system planning and resilience tasks, with a focus on capacity expansion planning problems. Furthermore, we plan to explore advanced methods for determining breakpoint locations and develop CPLAs for other quantities, such as current flows, as well as for systems with varying topologies.

\section*{Acknowledgment}
\noindent P. Buason thanks the Department of Mechatronics Engineering at Rajamangala University of Technology Phra Nakhon.

\IEEEtriggeratref{3}
\bibliographystyle{IEEEtran}
\bibliography{reference.bib}

\end{document}